\documentclass{article}

\usepackage{amsmath,amsfonts,amssymb}

\usepackage{graphicx}

\usepackage{dsfont}
\usepackage{bbold}
\usepackage{amsthm}
\usepackage[utf8]{inputenc}
\usepackage{appendix}

\newtheorem{thm}{Theorem}

\newtheorem{defn}[thm]{Definition}
\newtheorem{prp}{Proposition}

\newtheorem{conj}{Conjecture}

\title{\textbf{Sharp weak-type estimate for maximal operators associated to Cartesian families under an arithmetic condition}}

\author{Anthony Gauvan}

\begin{document}

\maketitle

\begin{abstract}
Given a set of integers $A \subset \mathbb{Z}$, we consider the smallest family $\mathcal{B}_{A^{n-1}}$ invariant by translation which contains the rectangles $$ R_{\boldsymbol{a}} = I_{a_1} \times \dots \times I_{a_{n-1}} \times I_{-(a_1+\dots+a_{n-1})}$$ for any $\boldsymbol{a} = (a_1,\dots,a_{n-1}) \in A^{n-1}$ and where $I_k = [0,2^k]$ for $k$ integer. We prove that if the set $A$ contains arbitrary large arithmetic progression then the maximal operator $M_{\mathcal{B}_{A^{n-1}}}$ associated to the family $\mathcal{B}_{A^{n-1}}$ is sharply bounded from $L^1\left(1+ \log^+ L^1 \right)^{n-1}$ to $L^{1,\infty}$.
\end{abstract}

We work in the Euclidean space $\mathbb{R}^n$ and for $U \subset \mathbb{R}^n$ we denote by $\left|U\right|$ its $n$-dimensional Lebesgue measure. We let $\mathcal{I}$ be the family containing every bounded intervals of the real line so that the family $\mathcal{I}^n$ defined as $$\mathcal{I}^n := \mathcal{I} \times \dots \times \mathcal{I}$$ stands for the family containing every axis parallel rectangles in $\mathbb{R}^n$. Given an integer $k \in \mathbb{Z}$, we let $I_k$ be the dyadic interval $[0,2^k]$ and we use the following notations: \begin{itemize}
    \item for $\boldsymbol{\hat{a}}$ in $\mathbb{Z}^n$, we consider the dyadic rectangle $R_{\boldsymbol{\hat{a}}}$ defined as$$R_{\boldsymbol{\hat{a}}} := I_{a_1} \times \dots \times I_{a_n},$$
    \item for $\boldsymbol{a}$ in $\mathbb{Z}^{n-1}$, we consider the dyadic rectangle $R_{\boldsymbol{a}}$ defined as$$R_{\boldsymbol{a}} := I_{a_1} \times \dots \times I_{a_{n-1}} \times I_{-(a_1+\dots+a_{n-1})}.$$
\end{itemize}

\section{Introduction}

The \textit{strong maximal operator} $M_n$ associated to the family $\mathcal{I}_n$ is a fundamental operator in harmonic analysis: it is defined for $f \in L^\infty$ and $x \in \mathbb{R}^n$ as $$M_n f(x) := \sup_{x \in R \in \mathcal{I}^n} \frac{1}{\left|R\right|} \int_R |f|.$$ In their seminal work \cite{JMZ}, Jessen, Marcinkiewicz, and Zygmund proved the so-called \textit{strong maximal Theorem} which specifies the sharp boundedness property of $M_n$ near $L^1$.

\begin{thm}[Strong maximal Theorem]\label{STRONG}
For any function $f$ and any $t > 0$, the following estimate holds $$ \left|\left\{M_nf > t \right\}  \right| \lesssim_n \int_{\mathbb{R}^n}  \frac{|f|}{t} \left(1 + \log^+ \frac{|f|}{t} \right)^{n-1}.$$ 
\end{thm}

It is not difficult to see that Theorem \ref{STRONG} is optimal in the following sense:  given any convex increasing function $\phi:[0, \infty) \rightarrow [0, \infty)$ such that $$\lim_{\infty} \frac{\phi(x)}{x (1 + \log^+ x)^{n-1}} = 0,$$ the operator $M_n$ cannot be bounded from the Orlicz space $L^\phi$ to $L^{1,\infty}$ that is to say we cannot have for any $f$ and $t$ $$ \left|\left\{M_nf > t \right\}  \right| \lesssim_n \int_{\mathbb{R}^n}  \phi\left(\frac{|f|}{t}\right).$$ Indeed, it is simple to check that given any integer $m \geq 1$, the following estimate holds $$\left| \left\{ M_n \mathbb{1}_Q \geq 2^{-m} \right\} \right| \gtrsim_n m^{n-1} 2^m \left|Q \right| $$ where $Q$ is any axis parallel cube in $\mathbb{R}^n$. Obviously more general maximal operators associated to geometric sets can be considered: given a family $\mathcal{B}$ included in $\mathcal{I}^n$ and invariant by translations, one can form in the same fashion the maximal operator $M_\mathcal{B}$ defined as $$M_\mathcal{B} f(x) := \sup_{x \in R \in \mathcal{B}} \frac{1}{\left|R\right|} \int_R |f|.$$ For example in $\mathbb{R}^3$, motivated by problematics arising in the theory of singular integrals,  Zygmund considered  the family $\mathcal{B}_Z$ composed of axis every axis parallel rectangles whose side lengths are of the form $$t \times s \times \sqrt{ts}$$ for any $t,s > 0$. Zygmund asked if one could prove that the maximal operator $M_{\mathcal{B}_Z}$ associated to the family $\mathcal{B}_Z$ is bounded from $L^1 \left(1 +\log^+ L^1 \right)$ to $L^{1,\infty}$ ? Loosely speaking, Zygmund expected such result because a rectangle in $\mathcal{B}_Z$, as a plane rectangle, is only defined by two parameters (up to translations). Cordoba answered positively to the question and the following conjecture was formulated.

\begin{conj}[Zygmund's conjecture I]
Let $\left\{ \phi_i : i \leq n\right\}$ be $n$ positive real functions depending on $k$ variables, increasing in each variable separately and assuming arbitrarily small values. Consider then basis $B_\phi$ of all axis parallel rectangles in $\mathbb{R}^n$ whose side lengths are of the form $$\phi_1(t) \times \dots \times \phi_n(t)$$ for any $t \in \mathbb{R}^k$. In this situation, the maximal operator $M_{B_{\boldsymbol{\phi}}}$ is bounded from $L^1 \left(1 +\log^+ L^1 \right)^{k-1}$ to $L^{1,\infty}$ 
\end{conj}

It appears that stated in this form, Zygmund's conjecture is false and this was proved by Soria in the simplest case $n=3$ and $k=2$. More recently, Rey exhibited new class of counter examples to this conjecture in \cite{R}. However, those negative results do not necessarily indicate that the idea behind Zygmund's conjecture is false but rather that it is not correctly formulated. Indeed, in his beautiful article \cite{S}, Stokolos proved the following Theorem (for the sake of clarity, we omit the \textit{geometric description} of this Theorem) thanks to the key idea of \textit{crystallisation} that he developed.

\begin{thm}[Stokolos]\label{TstoM}
Given any family $\mathcal{B}$ of axis parallel rectangles in $\mathbb{R}^2$ which is invariant by translations, there exists an integer $k \in \{1,2\}$ such that the maximal operator $M_\mathcal{B}$ is sharply bounded from $L^1\left(1 + \log^+ L^1 \right)^{k-1}$ to $L^{1,\infty}$.
\end{thm}

Here naturally, given an arbitrary positive number $k \geq 0$, we say that the maximal operator $M_\mathcal{B}$ is \textit{sharply bounded from $L^1\left( 1 \log^+ L^1\right)^k$ to $L^{1,\infty}$} when: \begin{itemize}
    \item For any function $f$ and any $t > 0$, the following estimate holds $$ \left|\left\{M_\mathcal{B}f > t \right\}  \right| \lesssim_{n,\mathcal{B}} \int_{\mathbb{R}^n}  \frac{|f|}{t} \left(1 + \log^+ \frac{|f|}{t} \right)^k.$$ 
    \item given any convex increasing function $\phi:[0, \infty) \rightarrow [0, \infty)$ such that $$\lim_{\infty} \frac{\phi(x)}{x (1 + \log^+ x)^k} = 0,$$ the operator $M_\mathcal{B}$ is not bounded from the Orlicz space $L^\phi$ to $L^{1,\infty}$.
\end{itemize} In regards of Theorem \ref{TstoM}, Stokolos proposed the following reformulation of Zygmund's conjecture in \cite{S2}.

\begin{conj}[Zygmund's conjecture II]
Given any family $\mathcal{B}$ invariant by translations included in $\mathcal{I}^n$, there exists an integer $1 \leq k \leq n$ such that the maximal operator $M_{\mathcal{B}}$ is sharply bounded from $L^1\left(1 + \log^+ L^1 \right)^{k-1}$ to $L^{1,\infty}$.
\end{conj}

So far, this conjecture has not been refuted and let us give non trivial examples: \begin{itemize}
    \item In $\mathbb{R}^n$, given an integer $1 \leq k \leq n$, consider the family $\mathcal{I}^n(k)$ defined as $$ \mathcal{I}^n(k) := \left\{ I_1 \times \dots \times I_n \in \mathcal{I}^n : \# \left\{ \left|I_j\right| :1 \leq j \leq n \right\} \leq k \right\}.$$ In \cite{Z}, Zygmund proved that the maximal operator $M_{n,k}$ associated to the family $\mathcal{I}^n(k)$ is bounded from $L^1\left(1+ \log^+ L^1 \right)^{k-1}$ to $L^{1,\infty}$ and easy computations show that those bounds are actually sharp.

    \item In $\mathbb{R}^n$, fix $n$ arbitrary infinite sets of integers $\left\{ A_i : 1 \leq i \leq n \right\}$ and denote by $\mathcal{B}$ the smallest family invariant by translations which contains the dyadic rectangles $$R_{\boldsymbol{a}} = I_{a_1} \times \dots \times I_{a_n}$$ for any $\boldsymbol{a} = (a_1,\dots,a_n) \in A_1 \times \dots \times A_n$. In \cite{S3}, Stokolos proved the maximal operator $M_{\mathcal{B}}$ associated to the family $\mathcal{B}$ is sharply bounded from $L^1\left(1 + \log^+ L^1 \right)^{n-1}$ to $L^{1,\infty}$. As explained by Stokolos, this classical result illustrates the fact that a simple \textit{rarefaction} of the family $\mathcal{I}^n$ into a family of the form $\mathcal{B}$ does not improve the boundedness property of the maximal operator $M_{\mathcal{B}}$ associated.
    
    \item In $\mathbb{R}^3$, given an infinite set of integers $S \subset \mathbb{N}$, consider the family $\mathcal{B}$ composed of every axis parallel rectangles whose side lengths are of the form $$ t \times s \times 2^j s$$ for any $t,s > 0$ and $j \in S$. In \cite{HS}, Hagelstein and Stokolos proved that the maximal operator $M_\mathcal{B}$ is sharply bounded sharply bounded from $L^1\left(1+ \log^+ L^1 \right)^2$ to $L^{1,\infty}$.

    \item In $\mathbb{R}^3$, given an infinite set of integers $S \subset \mathbb{N}$, consider the family $\mathcal{B}$ composed of every axis parallel rectangles whose side lengths are of the form $$ s \times \frac{2^j}{s} \times t$$ for any $t,s > 0$ and $j \in S$. In \cite{DHS}, Dmitrishin, Hagelstein and Stokolos proved that the maximal operator $M_\mathcal{B}$ is sharply bounded sharply bounded from $L^1\left(1+ \log^+ L^1 \right)^2$ to $L^{1,\infty}$.

\end{itemize}

\section{Result}

When studying maximal operators of the above type, a classic argument shows that given an arbitrary family $\mathcal{B}$ included in $\mathcal{I}^n$ and invariant by translations, one can suppose that any rectangle $R \in \mathcal{B}$ is - up to  translations - of the form $$ R_{\boldsymbol{\hat{a}}} = I_{a_1} \times \dots \times I_{a_n},$$ for some $\boldsymbol{\hat{a}} \in \mathbb{Z}^n$. Hence, given an arbitrary infinite set of integer $A \subset \mathbb{Z}$, let us denote by $\mathcal{B}_{A^{n-1}}$ the smallest family invariant by translations which contains the dyadic rectangle $$R_{\boldsymbol{a}} = I_{a_1} \times \dots \times I_{a_{n-1}} \times I_{-(a_1+\dots+a_{n-1})}$$ for any $\boldsymbol{a} = (a_1,\dots,a_{n-1}) \in A^{n-1}$. We will detail why it seems interesting to study families of the form $\mathcal{B}_{A^{n-1}}$ for infinite $A \subset \mathbb{Z}$ but observe for the moment that for any $\boldsymbol{a} \in A^{n-1}$, we have $$\left|R_{\boldsymbol{a}}\right| = 1.$$ We prove the following Theorem.

\begin{thm}\label{Tari}
If the set $A \subset \mathbb{Z}$ contains arbitrarily large arithmetic progressions then for any integer $m \geq 1$, there exists a set $E$ such that $$ \left| \left\{ M_{\mathcal{B}_{A^{n-1}}}\mathbb{1}_E \geq 2^{-m} \right\} \right| \gtrsim_n m^{n-1}2^m\left|E \right|.$$ In particular, the maximal operator $M_{\mathcal{B}_{A^{n-1}}}$ associated to the family $\mathcal{B}_{A^{n-1}}$ is sharply bounded from  $L^1\left(1 + \log^+ L^1 \right)^{n-1}$ to $L^{1,\infty}$.
\end{thm}

Thanks to the inclusion $\mathcal{B}_{A^{n-1}} \subset \mathcal{I}^n$ and the strong maximal Theorem \ref{STRONG}, it suffices to focus on the \textit{sharpness} of this bound. Let us discuss about the arithmetic condition: given any integer $m$, we suppose that $A$ contains a set $\left\{ u_0 < \dots < u_{m-1} \right\}$ such that for any $0 \leq k \leq m-1$ $$u_k = u_0 + k(u_1-u_0).$$ \textit{Let us insist that without this arithmetic condition on the set $A$, it is for the moment difficult to study the sharp boundedness property of the operator $M_{\mathcal{B}_{A^{n-1}}}$.} We will exploit the additive structure of $A$ in order to create a resonance on $n$-th axis: given an arbitrary $m \geq 1$, the set $E$ will be constructed as the Cartesian product of one dimensional \textit{crystals} \textit{i.e.} we will set $$E = C_1 \times \dots \times C_n.$$ A one dimensional crystal is a Cantor-like set whose structure is adapted to the specific scales ; we precise this in the following. If it is clear that the first crystals $\left\{ C_i : i \leq n-1 \right\}$ should be defined according to the set $A$, the last crystal $C_n$ must be defined with great care since we do not have freedom on this axis. Finally, the fact that $A$ contains arbitrarily large arithmetic condition is assured if, for example, it has a strictly positive asymptotic upper density thanks to Szemerédi's Theorem.

\section{Cartesian families}

The family $\mathcal{B}_{A^{n-1}}$ defined above is a specific type of \textit{Cartesian families} of rectangles and it would be desirable to obtain more informations on maximal operators associated to such families in order to progress on Zygmund's problem. More precisely, one could try to tackle Zygmund's problem with the \textit{additional hypothesis} that the family $\mathcal{B}$ is invariant by \textit{central dilations} \textit{i.e.} if for any $R \in \mathcal{B}$ and any $\lambda > 0$, we suppose that we have $$\lambda R \in \mathcal{B}.$$ In the same spirit than the dyadic reduction made above, one can show that such a family $\mathcal{B}$ is always generated by a family of dyadic rectangles of unit volume $$\left\{ R_{\boldsymbol{\hat{a}}} : \boldsymbol{\hat{a}} \in F \right\}$$ and so where the set $F \subset \mathbb{Z}^n$ is included in $$ F \subset \left\{ x_1 + \dots + x_n = 0 \right\}.$$ In the following, let us denote by $ D^{n,1}(0)$ the family $\left\{ x_1 + \dots + x_n = 0 \right\}$ that we identify with dyadic rectangles of unit volume anchored at $0$. One can easily inject $\mathbb{Z}^{n-1}$ into $D^{n,1}(0)$ as follow: given $\boldsymbol{a} \in \mathbb{Z}^{n-1}$, it suffices to define the rectangle $R_{\boldsymbol{a}} \in D^{n,1}(0)$ defined as  $$ R_{\boldsymbol{a}} := I_{a_1} \times \dots \times I_{a_{n-1}} \times I_{-(a_1+\dots+a_{n-1})}.$$ The fact that we have normalized along the $n$-th axis is not canonical but this is irrelevant for our problem. It would be interesting to give the sharp boundedness properties of maximal operators associated to \textit{Cartesian families}.

\begin{defn}[Cartesian family]
Given $(n-1)$ arbitrary sets of integers $$\left\{ A_i \subset \mathbb{Z} : 1 \leq i \leq n-1 \right\},$$ we say that the smallest family invariant by translations and central dilations that contains the family of dyadic rectangles of unit volume $$  \left\{ R_{\boldsymbol{a}} : \forall \boldsymbol{a} \in A_1 \times \dots \times A_{n-1} \right\} \subset D^{n,1}(0)$$ is a \textit{Cartesian family}. We denote this family $\mathcal{B}_{A_1 \times \dots \times A_{n-1}}$ and the maximal operator associated by $M_{\mathcal{B}_{A_1 \times \dots \times A_{n-1}}}$ or $M_{A_1 \times \dots \times A_{n-1}}$.
\end{defn}

It turns out that even in dimension $n =3$, it is difficult to handle Cartesian families. Shall we understand if Zygmund's problem is purely geometric or requires arithmetic tools, the following problems should be addressed:\begin{itemize}
    \item Can one specify the sharp boundedness property of the operator $M_{L \times \mathbb{Z}}$ where the set $L$ is defined as $$L = \left\{ 2^k : k \geq 1\right\}.$$ The set $L$ is a typical example of set of integers which do not have additive structure.
    \item Can one specify the sharp boundedness property of the operator $M_{L \times L}$ associated to the family $L \times L$ ?
    \item Generally, can one specify the sharp boundedness property of the operator $M_{A \times B}$  where $A$ and $B$ are arbitrary (infinite) sets of integers?
\end{itemize}

\section{Crystal}

The concept of \textit{crystallisation} was developed by Stokolos in \cite{S} and is fundamental for our purpose: we introduce the notion of \textit{crystal} in dimension $1$ and $n$ and we detail then different geometric properties of those sets that will be useful in the following.

Given an integer $a \in \mathbb{Z}$, we denote by $I_a$ the dyadic interval $[0,2^a]$ and by $O_a$ the \textit{oscillation at scale $a$} which is a subset of $\mathbb{R}$ defined as $$O_a := \bigcup_{k \in \mathbb{Z}} (k2^{a+1} + I_a).$$ A one dimensional \textit{crystal} is defined as follow.

\begin{defn}[Crystal]
Fix a finite set $A$ in $\mathbb{Z}$ that we denote $$A = \{ a_1 < \dots < a_m \}.$$ We define the crystal $C(A)$ as $$C(A) := I_{a_m} \cap \bigcap_{1 \leq i < m} O_{a_i}.$$
\end{defn}

An $n$-dimensional crystal is simply by definition the Cartesian product of $n$ one dimensional crystals \textit{i.e.} a set $E$ of the form $$E := C_1 \times \dots \times C_n.$$ We will usually denote a $n$-dimensional crystal by the letter $E$ or $Y$ and one dimensional crystal by $C$. We define the notion of \textit{primitive} rectangle associated to a $n$-dimensional crystal.

\begin{defn}[Primitive rectangle]
Given an $n$-dimensional crystal $E$, there exists a biggest dyadic rectangle $R$ anchored at $0$ included in $E$: we say that it is the \textit{primitive} rectangle of associated to the crystal $E$.
\end{defn}

We have the following disjointness property and we omit its proof since it is well known in the literature, see \cite{S3} for example.

\begin{prp}[Disjointness property]\label{PDP}
Fix a finite number of $n$-dimensional crystals $\left\{  E_i : i \leq N\right\}$ and suppose that the primitive rectangles associated $$\left\{  R_i : i \leq N\right\}$$ are independent \textit{i.e.} for any $i \leq N$ we have $$ \left| R_i - \bigcup_{j \neq i} R_i \right| \geq c_n \left| R_i \right|.$$ In this situation, we have the following estimate$$\left| \bigcup_{i \leq N} E_i \right| \simeq_n \sum_{i \leq N} \left|E_i \right|.$$
\end{prp}

In order to exploit $n$-dimensional crystals with a given maximal operator, we need to detail a specific property of those sets which indicate that, in some sense, they are well distributed at specific scales. The following notation will be useful: given a set of integer $ A = \left\{ a_1 < \dots < a_m \right\}$ and $1 \leq i \leq m$, we denote by $A[i]$ the set $$A[i] := \left\{ a_i < \dots <a_m \right\}.$$ Consider an $n$-dimensional crystal $E$ defined by $n$ sets of integers $\left\{A_k : k \leq n \right\}$ \textit{i.e.} the crystal $$ E = C(A_1) \times \dots \times C(A_n).$$ Fix then a multi-indice $\boldsymbol{\hat{i}} \in \mathbb{N}^n$ and consider the crystal $$Y(\boldsymbol{\hat{i}}) := C(A_1[i_1]) \times \dots \times C(A_n[i_n]).$$ Of coursed we have supposed that $i_k \leq \#A_k$. It is clear that we have $E \subset Y(\boldsymbol{\hat{i}})$ and denote then by $m$ the integer satisfying $$\left| Y(\boldsymbol{\hat{i}}) \right| = 2^m  \left| Y(\boldsymbol{\hat{i}}) \cap E\right|.$$  We claim the following.

\begin{prp}[Homogeneity property]\label{PHP}
If we denote by $R$ the primitive rectangle associated to $Y(\boldsymbol{\hat{i}})$, we have $$Y(\boldsymbol{\hat{i}}) \subset \left\{ M_R \mathbb{1}_E \geq 2^{-m} \right\}.$$ Here, $M_R$ stands for the maximal operator associated to the smallest family invariant by translations and which contains the dyadic rectangle $R$.
\end{prp}

\section{Proof of Theorem \ref{Tari}}

With Propositions \ref{PDP} and \ref{PHP} at hands, we are ready to prove Theorem \ref{Tari}. Given an arbitrary large integer $m \gg 1$, we construct an $n$-dimensional crystal $E \subset \mathbb{R}^{n}$ such that $$ \left| \left\{ M_{\mathcal{B}_{A^{n-1}}} \mathbb{1}_E > 2^{-m} \right\} \right| \gtrsim_n m^{n-1}2^m \left| E \right|.$$ By hypothesis, $A$ contains an arithmetic progression of length $m$ that we denote $$\left\{ u_0 < \dots < u_{m-1} \right\}.$$ For $s \in \{ 0, 1, \dots, m-1 \}$, define $$h_s := (n-1)u_0 + (u_1-u_0)s$$ and consider then the one dimensional crystals $X := C(u_0 < u_1 < \dots < u_{m-1})$ and $Z = C( -h_{m-1} < \dots < -h_0)$. Thanks to the crystal $X$ and $Z$, we define the $n$-dimensional crystal $E$ as $$E := X^{n-1}\times Z \subset \mathbb{R}^n.$$ Let us prove as claimed that we have $$\left| \left\{ M_{\mathcal{B}_{A^{n-1}}} \mathbb{1}_E > 2^{-m} \right\} \right| \gtrsim m^{n-1}2^m \left| E \right|.$$ Given a positive multi-indice $\boldsymbol{i} = (i_1, \dots ,i_{n-1})$ such that $$i_1 + \dots +i_{n-1} = s \leq m-1$$ consider the crystal $Y(\boldsymbol{i})$ defined as $$Y(\boldsymbol{i}) = C(u_{i_1} < \dots < u_{m-1}) \times  \dots \times C(u_{i_{n-1}} < \dots < u_{m-1}) \times C( -h_s < \dots < -h_0).$$ First observe that we have $E \subset Y(\boldsymbol{i})$ and also $$\left| Y(\boldsymbol{i}) \right| \simeq 2^{-i_1}\dots 2^{-i_{n-1}}  2^{-(m-1-s)}\left|E \right|  \simeq 2^{-m} \left|E \right|.$$ Observe now that the primitive rectangle associated to $Y(\boldsymbol{i})$ is the rectangle $$ R(\boldsymbol{i}) = I_{u_{i_1}} \times \dots \times I_{u_{i_{n-1}}} \times I_{-h_s}.$$ The crux of the argument lies in the fact that we have $$h_s = u_{i_1} + \dots + u_{i_{n-1}}$$ since $i_1 + \dots +i_{n-1} = s$ and that $\left\{ u_k : 1 \leq k \leq m \right\}$ is an arithmetic progression. Hence the primitive rectangle $R(\boldsymbol{i})$ satisfies $$ R(\boldsymbol{i}) \in \mathcal{B}_{A^{n-1}}$$ We apply now Proposition \ref{PHP} which yields $$Y(\boldsymbol{i}) \subset \left\{ M_{ R(\boldsymbol{i})} \mathbb{1}_E > 2^{-m} \right\} \subset \left\{ M_{\mathcal{B}_{A^{n-1}}} \mathbb{1}_E > 2^{-m} \right\}.$$ To conclude, it is not difficult to see that the family of rectangles $$\left\{  R(\boldsymbol{i}) : \boldsymbol{i} \geq 0, i_1 + \dots +i_{n-1} \leq m-1  \right\}$$ is independent and so applying Proposition \ref{PDP} we obtain $$\left| \bigcup_{\boldsymbol{i}} Y(\boldsymbol{i}) \right| \simeq \sum_{\boldsymbol{i}} \left|Y(\boldsymbol{i}) \right| \simeq m^{n-1} 2^m \left|E \right|.$$ This concludes the proof of Theorem \ref{Tari}.

{}

\end{document}